\documentclass[reqno,a4paper]{amsart}

\usepackage{amsmath,amsfonts,amssymb}
\usepackage{verbatim}
\usepackage{enumerate}
\usepackage{url}
\usepackage{tikz}
\usepackage{setspace}
\usepackage[left=2.5cm,right=2.5cm, top=3cm, bottom=2.5cm]{geometry}
\usepackage[utf8]{inputenc} 
\usepackage[T1]{fontenc}
\usepackage{hyperref} 	
\usetikzlibrary{arrows}

\def\d{{\sf d}}

\def\s{{\sf s}}
\def\e{{\sf \eta}}
\def\E{{\sf E}}
\def\F{{\mathcal F}}
\def\bd{{\boldsymbol{\cdot}}}

\def\ord{\mathop{\rm ord}\nolimits}

\def\lcm{\mathop{\rm lcm}}

\theoremstyle{plain}
\newtheorem{theorem}{Theorem}[section]
\newtheorem{lemma}[theorem]{Lemma}

\newtheorem{proposition}[theorem]{Proposition}

\newtheorem{conjecture}[theorem]{Conjecture}
\def\proof{\noindent {\it Proof: }}
\def\qed{\hfill\hbox{$\square$}}

\theoremstyle{definition}

\numberwithin{equation}{section}

\title{The main zero-sum constants over $D_{2n} \times C_2$}

\author[F. E. Brochero Mart\'{\i}nez]{F. E. Brochero Mart\'{\i}nez}
\address{
Departamento de Matem\'{a}tica\\
Universidade Federal de Minas Gerais (UFMG)\\
Belo Horizonte, MG\\
31270-901\\
Brazil\\
}
\email{fbrocher@mat.ufmg.br }

\author[A. Lemos]{A. Lemos}
\address{Departamento de Matem\'{a}tica\\
Universidade Federal de Vi\c cosa (UFV)\\
Vi\c cosa, MG\\
36570-000\\
Brazil\\
}
\email{abiliolemos@ufv.br}

\author[B. K. Moriya]{B. K. Moriya}
\address{Departamento de Matem\'{a}tica\\
Universidade Federal de Vi\c cosa (UFV)\\
Viçosa, MG\\
36570-000\\
Brazil\\
}
\email{bhavinkumar@ufv.br}

\author[S. Ribas]{S. Ribas}
\address{
Departamento de Matem\'{a}tica\\
Universidade Federal de Ouro Preto (UFOP)\\
Ouro Preto, MG\\
35400-000\\
Brazil\\
}
\email{savio.ribas@ufop.edu.br }

\thanks{The first named author was partially supported by FAPEMIG APQ-02973-17, Brazil.}

\date{\today}

\onehalfspace

\usepackage{xcolor}

\begin{document}

\maketitle

\begin{abstract}
Let $C_2$ be the cyclic group of order $2$ and $D_{2n}$ be the dihedral group of order $2n$, where $n$ is even. In this paper, we provide the exact values of some zero-sum constants over $D_{2n} \times C_2$, namely small Davenport constant, Gao constant, $\e$-constant and Erd\H os-Ginzburg-Ziv constant. As a consequence, we prove the Gao's and Zhuang-Gao's Conjectures for this group. These are the first concrete results on zero-sum problems for a family of non-abelian groups of rank greater than $2$.
\end{abstract}

\section{Introduction}

Let $G$ be a finite multiplicative group. The {\em zero-sum problems} consist of studying the conditions which ensures that a given sequence over $G$ has a non-empty product-one subsequence with some prescribed properties. This kind of problem dates back to the pioneering works of Erd\H os, Ginzburg \& Ziv \cite{EGZ}, van Emde Boas \& Kruyswijk \cite{vEBK} and Olson \cite{Ols1,Ols2}. It has applications and connections in several branches of mathematics, such as number theory \cite{AGP}, coding theory \cite{PS}, factorization theory \cite{CFGO}, and finite geometry \cite{Ed,EFLS}. For an overview on zero-sum problems over abelian groups one can refer to the surveys from Caro \cite{Car} and Gao \& Geroldinger \cite{GaGe}. In the 80s, the zero-sum problems were further generalized to non-abelian groups (see, for instance, \cite{Yu,YP,ZhGa}), and this explains the use of multiplicative notation; in particular, the use of product-one rather than zero-sum.

\subsection{Definitions and notations}

Let $\F(G)$ be a free abelian monoid, written multiplicatively, with basis $G$. A sequence $S$ over a finite group $G$ is an element of $\F(G)$. Note that $\F(G)$ is equipped with the sequence concatenation product denoted by ${\bd}$. A sequence $S \in \F(G)$ has the form 
$$S = g_1 {\bd} \dots {\bd} g_k = \prod_{1\le i \le k} {\color{white}.}\!\!\!\!\!\!^{\bullet} \, g_i = \prod_{1\le i \le k} {\color{white}.}\!\!\!\!\!\!^{\bullet} \,  g_{\tau(i)},$$ 
for any permutation $\tau: \{1,2,\dots,k\} \to \{1,2,\dots,k\}$, where $g_1, \dots, g_k \in G$ are the {\em terms} of $S$ and $k = |S| \ge 0$ is the {\em length} of $S$. Given $g \in G$ and $t \ge 0$, we abbreviate $g^{[t]} = \prod_{1\le i \le t}^{\bullet} g$. For $g \in G$, the {\em multiplicity} of the term $g$ in $S$ is denoted by $v_g(S) = \#\{i \in \{1,2,\dots,k\} \mid g_i = g\}$, therefore we may also write $S = \prod_{g \in G}^{\bullet} g^{[v_g(S)]}$. A sequence $T$ is a {\em subsequence} of $S$ if $v_g(T) \le v_g(S)$ for all $g \in G$. In this case, we use the notation $T \mid S$ and write $S {\bd} T^{-1} = \prod_{g \in G}^{\bullet} g^{[v_g(S) - v_g(T)]}$ and $T^{[k]} = \prod_{1\le i \le k}^\bullet T$. Moreover, let:
\begin{align*}
\pi(S) &= \{g_{\tau(1)} \dots g_{\tau(k)} \in G \mid \tau \text{ is a permutation of $\{1,2,\dots,k\}$}\} \;\;\; \text{ be the {\em set of products} of $S$}; \\
\Pi(S) &= \bigcup_{T \mid S \atop |T| \ge 1} \pi(T) \subset G \;\;\; \text{ be the {\em set of subsequence products} of $S$}; \\
S \cap K &= \prod_{g \mid S \atop g \in K} {\color{white}.}\!\!\!\!\!^{\bullet} \, g \;\;\; \text{ be the subsequence of $S$ formed by the terms that lie in a subset $K$ of $G$}. 
\end{align*}
\newpage

The sequence $S$ is called:
\begin{itemize}
\item {\em product-one free} if $1 \not\in \Pi(S)$,
\item {\em product-one sequence} if $1 \in \pi(S)$,
\item {\em minimal product-one sequence} if it is a non-empty product-one sequence and every proper subsequence is product-one free,
\item {\em short product-one sequence} if it is a product-one sequence of length $1 \le |S| \le \exp(G) = \lcm\{\ord(g): g \in G\}$.
\end{itemize}

For a given finite multiplicative group $G$, the following zero-sum invariants are defined:
\begin{itemize}
\item {\em Small Davenport constant}, $\d(G)$, is the maximum length of a product-one free sequence over $G$,
\item {\em Gao constant}, $\E(G)$, is the smallest $\ell > 0$ such that every $S \in \F(G)$ with $|S| \ge \ell$ has a product-one subsequence of length $|G|$,
\item {\em $\e$-constant}, $\e(G)$, is the smallest $\ell > 0$ such that every $S \in \F(G)$ with $|S| \ge \ell$ has a short product-one subsequence,
\item {\em Erd\H os-Ginzburg-Ziv constant}, $\s(G)$, is the smallest $\ell > 0$ such that every $S \in \F(G)$ with $|S| \ge \ell$ has a product-one subsequence of length $\exp(G)$.
\end{itemize}
It is possible to show that they are all well-defined and finite for finite groups $G$.

\subsection{Background on the zero-sum constants}

Let $C_n$ be the cyclic group of order $n$ and $D_{2n}$ be the dihedral group of order $2n$. The following table summarizes some known results with references (not necessarily the main or the first):

\vspace{2mm}
\begin{center}
\begin{tabular}{|c|c|c|c|c|}
\hline
$G$ & $\d(G)$ & $\E(G)$ & $\e(G)$ & $\s(G)$ \\
\hline
\hline
$C_n$ & $n-1$ & $2n-1$ \cite{EGZ} & $n$ & $2n-1$ \cite{EGZ} \\
\hline
$C_m \times C_n$, $m \mid n$ & $m+n-2$ \cite{Ols2} & $\d(G)+|G|$ \cite{Ga} & $2m+n-2$ \cite{GeHK} & $2m+2n-3$ \cite{GeHK} \\
\hline
$C_{p^{\alpha_1}} \times \dots \times C_{p^{\alpha_k}}$, & $\sum_{i=1}^k (p^{\alpha_i} - 1)$ & $\d(G) + |G|$ \cite{Ga} & unknown & unknown \\
$p$ prime, $\alpha_1 \le \dots \le \alpha_k$ & \cite{Ols1} &  &  &  \\
\hline
$D_{2n}$, $n$ even & $n$ \cite{Bas} & $3n$ \cite{Bas} & $n+1$ \cite{Zh} & $2n$ \cite{OhZh} \\
$D_{2n}$, $n$ odd & (both parities) & (both parities) & (both parities) 
& $3n$ \cite{OhZh} \\
\hline
\end{tabular}
\end{center}
\vspace{2mm}

The scenario is still obscure for groups of rank greater than two (see, for example, \cite{Liu} and the references therein, where the results broke the expected rule). Concrete results are known only for few families of groups (see, for instance, \cite{Sc1} and references therein).

By Pigeonhole Principle, it holds $\d(G) +1 \le |G|$ for every group $G$. Olson \& White \cite{OW} proved that if $G$ is non-abelian then $\d(G) \le \lfloor |G|/2 \rfloor$. We also have that $\d(G) \le \e(G)$ by definition, and in addition $\E(G) \le 2|G|-1$ (see \cite[Theorem~10.1]{Gr1}). Some routine arguments (see \cite[Lemma 4]{ZhGa}) yield 
\begin{equation}\label{sexp}
\s(G) \ge \e(G)+\exp(G)-1
\end{equation}
and
\begin{equation}\label{Ed}
\E(G) \ge \d(G) + |G|.
\end{equation}

\begin{conjecture}[Gao \cite{Ga2}]\label{conjgao}
$\s(G) = \e(G) + \exp(G) - 1$ for every finite group $G$.
\end{conjecture}

\begin{conjecture}[Zhuang-Gao \cite{ZhGa}]\label{conj1}
$\E(G) = \d(G) + |G|$ for every finite group $G$.
\end{conjecture}

\subsection{On the group $D_{2n} \times C_2$}

This paper deals with some zero-sum constants over the non-abelian group of rank three:
$$D_{2n} \times C_2 = \langle x,y,z \mid x^2 = y^n = z^2 = 1, \; yx = xy^{-1}, \; xz = zx, \; yz = zy \rangle,$$ 
where $n \ge 4$ is even. We observe more generally that if $\gcd(m,n)=1$, then, by \cite[Section~1.8]{Hung}, $D_{2n} \times C_m \simeq D_{2mn}$. Hence, by the Fundamental Theorem of Finitely Generated Abelian Groups, we may assume without loss of generality that $m \mid n$ and that $m > 1$. It explains why we set $n$ even for $m=2$.

Let $m \ge 2$ and $n \ge 3$ be integers such that $m \mid n$. We expect that 
\begin{align*}
\d(D_{2n} \times C_m) &= m + n - 1, \\
\E(D_{2n} \times C_m) &= 2mn + m + n - 1, \\
\e(D_{2n} \times C_m) &= 2m + n - 1, \\
\s(D_{2n} \times C_m) &= 2m + 2n - 2.
\end{align*}

Some simple arguments show that the lower bounds of in the previous expressions hold (see Lemma \ref{lowerbound} and Ineqs. \eqref{sexp} and \eqref{Ed}). In the particular case that $m=2$, we are able to prove the respective upper bounds. It is worth mentioning that these are the first concrete results on the zero-sum problems for non-abelian groups of rank greater than two.

The paper is organized as follows. In Section \ref{lemmas}, we present some auxiliary results that will be used throughout the paper, which include general lower bounds for $\d(G)$ and $\e(G)$, and a result on zero-sum problems. In Sections \ref{conjseexp} and \ref{conjEdord}, respectively, we prove Conjectures \ref{conjgao} and \ref{conj1} for the group $G \simeq D_{2n} \times C_2$.

\section{The lower bounds and another auxiliary result}\label{lemmas}

In this section, we present the auxiliary results that will be used throughout the paper. We start with the following:

\begin{lemma}\label{lemineq}
Let $G$ be a finite group and let $H$ be a normal subgroup of $G$. We have that 
\begin{enumerate}[(a)]
\item $$\d(G) \ge \d(H) + \d(G/H).$$
\item If $\exp(H) = \exp(G)$, then $$\e(G) \ge \e(H) + \d(G/H).$$ 
\end{enumerate}
\end{lemma}

\proof
\begin{enumerate}[{\it(a)}]
\item Let $S_1 \in \F(G)$ with $|S_1| = \d(G/H)$ be a sequence without subsequences whose products (in any order) belong to $H$ (that is, $S_1$ is product-one free over the quotient $G/H$), and let $S_2 \in \F(H)$ with $|S_2| = \d(H)$ be a product-one free sequence. It follows that $S = S_1 \bd S_2 \in \F(G)$ is product-one free and $|S| = \d(H) + \d(G/H)$.
\item Let $S_1 \in \F(G)$ be as in the previous case, and let $S_2 \in \F(H)$ be a sequence with no short product-one subsequences such that $|S_2| = \e(H) - 1$. By ``short'' here we mean sequences of length at most $\exp(G) = \exp(H)$. It follows that $S = S_1 \bd S_2 \in \F(G)$ satisfies $|S| = \d(G/H) + \e(H) - 1$ and we claim that $S$ has no short product-one subsequences. In fact, suppose that $T = T_1 \bd T_2 \mid S$ is a short product-one sequence, where $T_1 \mid S_1$ and $T_2 \mid S_2$. Since $\pi(T) = 1$, we have that $\pi(T_1) \in H$, therefore $T_1$ is empty. Hence $T = T_2$, but this is a contradiction by definition of $T$.
\end{enumerate}

\qed

\vspace{2mm}

We now present the lower bounds that have motivated our expectations on the small Davenport constant and on the $\e$-constant.

\begin{lemma}\label{lowerbound}
Let $m \ge 2$ and $n \ge 3$ be integers such that $m \mid n$. We have that
\begin{enumerate}[(a)]
\item $$\d(D_{2n} \times C_m) \ge m+n-1.$$
\item If $n$ is even, then $$\e(D_{2n} \times C_m) \ge 2m+n-1.$$
\end{enumerate}
\end{lemma}

\proof
Let $G \simeq D_{2n} \times C_m$ and $H \simeq C_n \times C_m$, so that $\exp(H) = \exp(G)$ provided $n$ is even. Since the index $[G:H] = 2$, it follows that $H$ is a normal subgroup of $G$. This lemma now follows from the previous lemma.

\qed

\vspace{2mm}

Sometimes, we need to ensure that a given sequence over $G$ of large enough length is not product-one free. The study of the product-one free sequences of large length (with some prescribed property) is called {\em inverse zero-sum problem}. The next result is concerned with the inverse problem related to Erd\H os-Ginzburg-Ziv constant over a cyclic group.

\begin{lemma}[\cite{EGZ,Ga1}]\label{inverseegz}
For every integer $n \ge 2$, $\s(C_n) = 2n-1$. In addition, let $2 \le k \le \lfloor n/2 \rfloor + 2$, and suppose that $S \in \F(C_n)$ and $|S| = 2n-k$. If every subsequence $T \mid S$ of length $|T| = n$ is product-one free, then there exist distinct elements $a \mid S$ and $b \mid S$ such that $\min\{v_a(S), v_b(S)\} \ge n - 2k + 3$ and $v_a(S) + v_b(S) \ge 2n - 2k + 2$, where $ab^{-1}$ generates $C_n$.
\end{lemma}

From now on, we let $n \ge 4$ be an even integer and $m=2$. We will also use the following notation:

\vspace{2mm}

\noindent {\bf Notations.} {\it For $S \in \F(D_{2n} \times C_2)$, we write $S = S_1 \bd S_2 \bd S_3 \bd S_4$, where $$S_1 = \prod_i {\color{white}.}\!\!\!^{\bullet} (xy^{u_i}z), \; S_2 = \prod_i {\color{white}.}\!\!\!^{\bullet} (xy^{v_i}), \; S_3 = \prod_i {\color{white}.}\!\!\!^{\bullet} (y^{w_i}z), \; \text{ and } S_4 = \prod_i {\color{white}.}\!\!\!^{\bullet} (y^{t_i}),$$
and $|S| = |S_1| + |S_2| + |S_3| + |S_4|$.

In addition, for $g = x^{\alpha}y^{\beta}z^{\gamma} \in D_{2n} \times C_2$, where $\alpha, \gamma \in \{0,1\}$ and $0 \le \beta \le n-1$, we define the {\em parity} of $g$ in the following way: The parity of $g$ is {\em even} if $\beta$ is even, and the parity of $g$ is {\em odd} if $\beta$ is odd.}

\vspace{3mm}

The following lemma helps us to solve the small cases of Propositions \ref{sle} and \ref{Ele}, namely, $n \in \{4,6,8\}$.

\begin{lemma}\label{termosR}
Let $n \ge 4$ be an even integer and let $S \in \F(D_{2n} \times C_2)$ such that $|S_i| \le 2$ for every $i \in \{1,2,3,4\}$, where the terms of each $S_i$ have opposite parities provided $|S_i| = 2$. If $5 \le |S| \le 8$, then there exists $T \mid S$ with $|T| = 4$ whose product belongs to $C_{n/2} \simeq \langle y^2 \rangle$. Moreover, 
\begin{enumerate}[{\em (a)}]
\item if $6 \le |S| \le 8$, then $\pi(T)$ has at least two distinct elements provided either $n/2$ is even or, for every integer $\alpha$, $S \nmid xy^{\alpha}z \bd xy^{\alpha+n/2}z \bd xy^{\alpha} \bd xy^{\alpha+n/2} \bd z \bd y^{n/2}z \bd 1 \bd y^{n/2}$;

\item if $|S| = 5$ and $|S_1| + |S_2| \ge 2$, then $\pi(T)$ has at least two distinct elements provided either $n/2$ is even or, for every integer $\alpha$, $T \nmid xy^{\alpha}z \bd xy^{\alpha+n/2}z \bd xy^{\alpha} \bd xy^{\alpha+n/2} \bd z \bd y^{n/2}z \bd 1 \bd y^{n/2}$.
\end{enumerate}
\end{lemma}

\proof
We split into two cases and some subcases.
\begin{enumerate}[(a)]
\item {\bf Case} $6 \le |S| \le 8$. Since $|S_i| \le 2$, we have that $2 \le |S_1| + |S_2| \le 4$. 
\begin{enumerate}[{(a.}1)]
\item {\bf Subcase} $|S_1| + |S_2| = 4$. Let
$$T = xy^{u_1}z \bd xy^{u_2}z \bd xy^{v_1} \bd xy^{v_2},$$ 
so that $y^{u_2 - u_1 + v_2 - v_1}, y^{u_2 - u_1 + v_1 - v_2} \in \pi(T)$. Since $u_1 \not\equiv u_2 \pmod 2$ and $v_1 \not\equiv v_2 \pmod 2$, these products belong to $C_{n/2}$. These products are further distinct if and only if $v_1 \not\equiv v_2 \pmod {n/2}$. Since $v_1 \not\equiv v_2 \pmod 2$, we are done unless $n/2$ is odd. In this case, changing the order of $xy^{u_1}z$ and $xy^{u_2}z$, we obtain similarly that either the products are distinct or $u_1 \equiv u_2 \pmod {n/2}$. Changing the order of $xy^{u_2}z$ and $xy^{v_1}$, we also obtain that either the products are distinct or (without loss of generality) $u_1 \equiv v_1 \equiv u_2 + n/2 \equiv v_2 + n/2 \pmod n$. Hereupon, we observe that $\pi(T) = \{1\}$.

If $|S_3| = 2$, say $y^{w_1}z \bd y^{w_2}z \mid S$, then let 
$$T' = xy^{\alpha}z \bd xy^{\alpha+n/2}z \bd y^{w_1}z \bd y^{w_2}z,$$ 
so that $\pi(T') \subset C_{n/2}$. Since $y^{\pm w_1 \pm w_2+n/2} \in \pi(T')$, we have $y^{\pm w_1 \pm w_2+n/2} = 1$ if and only if $w_1 \pm w_2 \equiv n/2 \pmod n$. It implies that $\{w_1, w_2\} = \{0, n/2\}$ modulo $n$. Similarly, we obtain $\{t_1, t_2\} = \{0, n/2\}$ modulo $n$ provided $|S_4| = 2$, that is, $y^{t_1} \bd y^{t_2} \mid S$.

If $|S_3| = |S_4| = 1$, say $y^{w_1}z \bd y^{t_1} \mid S$, then let 
$$T'' = 
\begin{cases}
xy^{\alpha}z \bd xy^{\alpha} \bd y^{w_1}z \bd y^{t_1} &\text{ if $w_1 \equiv t_1 \pmod 2$}, \\
xy^{\alpha}z \bd xy^{\alpha+n/2} \bd y^{w_1}z \bd y^{t_1} &\text{ if $w_1 \not\equiv t_1 \pmod 2$},
\end{cases}$$
so that $\pi(T'') \subset C_{n/2}$. If $w_1 \equiv t_1 \pmod 2$, then changing orders and imposing $\pi(T'') = 1$, we obtain $w_1 \pm t_1 \equiv 0 \pmod n$, which implies either $w_1\equiv t_1 \equiv 0 \pmod n$ or $w_1 \equiv t_1 \equiv n/2 \pmod n$. If $w_1 \not\equiv t_1 \pmod 2$, then $w_1 \pm t_1 \equiv n/2 \pmod n$, which implies $\{w_1, t_1\} = \{0, n/2\}$ modulo $n$.

\item {\bf Subcase} $|S_1| + |S_2| = 3$. Then either $|S_1| = 2$ or $|S_2| = 2$. These two cases are completely similar, therefore we just need to prove for $|S_1| = 2$. Since $|S| \ge 6$, we also have either $|S_3| = 2$ or $|S_4| = 2$. Without loss of generality, suppose that $|S_3| = 2$. Let 
$$T = xy^{u_1}z \bd xy^{u_2}z \bd y^{w_1}z \bd y^{w_2}z,$$ 
so that $y^{u_2 - u_1 \pm w_1 \pm w_2} \in \pi(T) \subset C_{n/2}$. Since $u_1$ and $u_2$ have opposite parities, as well as $w_1$ and $w_2$ do, these products belong to $C_{n/2}$. Moreover, these products are equal if and only if $n/2$ is odd, $u_1 \equiv u_2 + n/2 \pmod n$ and $\{w_1, w_2\} = \{0, n/2\}$ modulo $n$. In this case, let 
$$T' = 
\begin{cases}
xy^{u_1}z \bd xy^{v_1} \bd y^{w_1}z \bd y^{t_1} &\text{ if $u_1+v_1+w_1+t_1$ is even}, \\
xy^{u_2}z \bd xy^{v_1} \bd y^{w_1}z \bd y^{t_1} &\text{ if $u_2+v_1+w_1+t_1$ is even},
\end{cases}$$
so that $\pi(T') \subset C_{n/2}$. The same argument as previous subcase implies that $v_1 \in \{u_1,u_2\}$ and $t_1 \in \{0, n/2\}$ modulo $n$, therefore we are done.

\item {\bf Subcase} $|S_1| + |S_2| = 2$. Then either $(|S_1|,|S_2|) = (2,0)$, $(|S_1|,|S_2|) = (1,1)$, or $(|S_1|,|S_2|) = (0,2)$. The first and third cases are completely similar, therefore we just need to prove the first and the second. In any case, $|S_3| = |S_4| = 2$.

For $(|S_1|,|S_2|) = (2,0)$, we take 
$$T = 
\begin{cases}
xy^{u_1}z \bd xy^{u_2}z \bd y^{w_1}z \bd y^{w_2}z, \text{ or } \\
xy^{u_1}z \bd xy^{u_2}z \bd y^{t_1} \bd y^{t_2},
\end{cases}$$ 
so that $\pi(T) \subset C_{n/2}$. The same argument as Subcase (a.1) for $|S_3|=2$ works.

For $(|S_1|,|S_2|) = (1,1)$, we take either
$$T' = 
\begin{cases}
xy^{u_1}z \bd xy^{v_1} \bd y^{w_1}z \bd y^{t_1}, &\text{ if $u_1+v_1+w_1+t_1$ is even, or }\\
xy^{u_1}z \bd xy^{v_1} \bd y^{w_1}z \bd y^{t_2}, &\text{ if $u_1+v_1+w_1+t_2$ is even},
\end{cases}$$ 
so that $\pi(T') \subset C_{n/2}$. The same argument than Subcase (a.2) for its $T'$ works, therefore we are done.
\end{enumerate}

\item {\bf Case} $|S| = 5$. If $|S_1| + |S_2| = 1$, then the only possible $T$ is $T = y^{w_1}z \bd y^{w_2}z \bd y^{t_1} \bd y^{t_2}$, so that $\pi(T) = \{y^{w_1+w_2+t_1+t_2}\} \subset C_{n/2}$. On the other hand, suppose that $|S_1| + |S_2| \ge 2$. By Pigeonhole Principle, there exists $i \in \{1,2,3,4\}$ such that $|S_i| = 2$. We consider the following subcases:
\begin{enumerate}[{(b.}1)]
\item {\bf Subcase} $|S_1| = |S_2| = 2$. The only possible $T$ is
$$T = xy^{u_1}z \bd xy^{u_2}z \bd xy^{v_1} \bd xy^{v_2},$$
so that $\pi(T) \subset C_{n/2}$. By the same arguments of Case (a.1), we have that either $\pi(T)$ has two distinct elements, or $u_1 \equiv u_2 + n/2 \pmod n$ and $w_1 \equiv w_2 + n/2 \pmod n$. The fifth term of $S$ can be anyone in $S_3 \bd S_4$.

\item {\bf Subcase} $|S_1| = |S_3| = 2$. The only possible $T$ is
$$T = xy^{u_1}z \bd xy^{u_2}z \bd y^{w_1}z \bd y^{w_2}z,$$
so that $\pi(T) \subset C_{n/2}$. By the same arguments of Case (a), we have that either $\pi(T)$ has two distinct elements, or $u_1 \equiv u_2 + n/2 \pmod n$ and $w_1 \equiv w_2 + n/2 \pmod n$. The fifth term of $S$ can be anyone in $S_2 \bd S_4$.

The subcases $|S_1|=|S_4|=2$, $|S_2|=|S_3|=2$ and $|S_2|=|S_4|=2$ follow similarly.

\item {\bf Subcase} $|S_1| = 2$, $|S_2| = |S_3| = |S_4| = 1$. The only possible $T$ is
$$T = 
\begin{cases}
xy^{u_1}z \bd xy^{v_1} \bd y^{w_1}z \bd y^{t_1}, &\text{ if $u_1+v_1+w_1+t_1$ is even, or }\\
xy^{u_2}z \bd xy^{v_1} \bd y^{w_1}z \bd y^{t_1}, &\text{ if $u_2+v_1+w_1+t_1$ is even},
\end{cases}$$ 
so that $\pi(T) \subset C_{n/2}$. As the $T'$ of Subcase (a.2), either $\pi(T)$ has two distinct elements, or (say $u_1+v_1+w_1+t_1$ is even) $u_1 - v_1 \pm w_1 \pm t_1 \equiv 0 \pmod n$ for every choice of signs and $n/2$ is odd. The latter implies that either $u_1 \equiv v_1 \pmod n$ or $u_1 \equiv v_1 + n/2 \pmod n$, and further $w_1,t_1 \in \{0, n/2\}$ modulo $n$. The fifth term of $S$ is $xy^{u_2}z \mid S_1$, where $u_2 \not\equiv u_1 \pmod 2$.

The subcase $(|S_1|,|S_2|,|S_3|,|S_4|) = (1,2,1,1)$ follows similarly.

\item {\bf Subcase} $|S_3| = 2$, $|S_1| = |S_2| = |S_4| = 1$. The only possible $T$ is
$$T = 
\begin{cases}
xy^{u_1}z \bd xy^{v_1} \bd y^{w_1}z \bd y^{t_1}, &\text{ if $u_1+v_1+w_1+t_1$ is even, or }\\
xy^{u_1}z \bd xy^{v_1} \bd y^{w_2}z \bd y^{t_1}, &\text{ if $u_1+v_1+w_2+t_1$ is even},
\end{cases}$$ 
so that $\pi(T) \subset C_{n/2}$. As the $T'$ of Subcase (a.3), either $\pi(T)$ has two distinct elements, or (say $u_1+v_1+w_1+t_1$ is even) $u_1 - v_1 \pm w_1 \pm t_1 \equiv 0 \pmod n$ for every choice of signs and $n/2$ is odd. The latter implies that either $u_1 \equiv v_1 \pmod n$ or $u_1 \equiv v_1 + n/2 \pmod n$, and further $w_1,t_1 \in \{0, n/2\}$ modulo $n$. The fifth term of $S$ is $y^{w_2}z \mid S_3$, where $w_2 \not\equiv w_1 \pmod 2$.

The subcase $(|S_1|,|S_2|,|S_3|,|S_4|) = (1,1,1,2)$ follows similarly.
\end{enumerate}
\end{enumerate}
It completes the proof of the lemma.

\qed

\section{Proof of Conjecture $\s(G) = \e(G) + \exp(G) - 1$ for $G \simeq D_{2n} \times C_2$}\label{conjseexp}

\begin{proposition}\label{sle}
Let $n \ge 4$ be an even integer. Then $\s(D_{2n} \times C_2) \le 2n+2$.
\end{proposition}

\proof
Let $S \in \F(D_{2n} \times C_2)$ with $|S| = 2n+2$. If $|S_1| + |S_2| \le 1$, then $|S_3| + |S_4| \ge 2n+1 = \s(C_n \times C_2)$. Since $S_3 \bd S_4 \in \F(C_n \times C_2)$, we are done. From now on, we assume that $|S_1|+|S_2| \ge 2$. If $|S_i| \ge 3$ for some $i \in \{1,2,3,4\}$, then there exists $g_{i_1} \bd g_{i_2} \mid S_i$ with the same parity, which implies that $g_{i_1} \cdot g_{i_2} \in C_{n/2} = \langle y^2 \rangle$. This argument can be repeated $\left\lceil \frac{|S_i| - 2}{2} \right\rceil$ times. Therefore, we obtain at least
\begin{equation}\label{eqceiln-3}
\left\lceil \frac{|S_1| - 2}{2} \right\rceil + \left\lceil \frac{|S_2| - 2}{2} \right\rceil + \left\lceil \frac{|S_3| - 2}{2} \right\rceil + \left\lceil \frac{|S_4| - 2}{2} \right\rceil \ge \frac{|S|-8}{2} = n - 3 = 2 \cdot \frac{n}{2} - 3
\end{equation}
terms over $C_{n/2}$, each of them being a product of exactly two terms of $S_i$, for $i \in \{1,2,3,4\}$. Let $S' \in \F(C_{n/2})$ denote these (at least) $n-3$ terms. We look at the terms of $S$ that are not used to form $S'$, and denote this sequence by $R$. It follows that the terms of $R$ lying in the same $S_i$, for $i \in \{1,2,3,4\}$, have opposite parities, otherwise we would obtain one more pair whose product lies in $C_{n/2}$. We have $|R| \in \{0,2,4,6,8\}$ and the left hand side of Ineq. \eqref{eqceiln-3} equals $\frac{2n+2-|R|}{2} = n+1-\frac{|R|}{2}$. If $|R| \le 4$, then the left hand side of Ineq. \eqref{eqceiln-3} is at least $n-1 = \s(C_{n/2})$, therefore we are done. Thus we assume that $|R| \in \{6,8\}$. 

If $|R| = 6$, then the left hand side of Ineq. \eqref{eqceiln-3} can be improved to $n-2$. If $S'$ has no product-one subsequence of length $n/2$, then Lemma \ref{inverseegz} ensures that 
$$S' = (y^{2\alpha})^{[n/2-1]} \bd (y^{2\beta})^{[n/2-1]},$$ 
where $\gcd(\alpha-\beta,n/2) = 1$.

If $|R| = 8$, then the left hand side of Ineq. \eqref{eqceiln-3} equals $n-3$. If $S'$ has no product-one subsequence of length $n/2$, then Lemma \ref{inverseegz} ensures that either
\begin{align*}
S' &= (y^{2\alpha})^{[n/2-2]} \bd (y^{2\beta})^{[n/2-1]}, \\ 
S' &= (y^{2\alpha})^{[n/2-2]} \bd (y^{2\beta})^{[n/2-2]} \bd (y^{2\gamma}), \quad \text{or} \\ 
S' &= (y^{2\alpha})^{[n/2-3]} \bd (y^{2\beta})^{[n/2-1]} \bd (y^{2\gamma}),
\end{align*}
where $\gcd(\alpha-\beta, n/2) = 1$. Notice that if $n \le 6$, then the third sequence becomes $S' = (y^{2\beta})^{[n/2-1]} \bd (y^{2\gamma})$, where $\gcd(\gamma-\beta,n/2) = 1$.

We now use mainly the following strategy: select two or four terms involving at least two in $S_1 \bd S_2$ to yield an element in $C_{n/2}$, place other elements of $S'$ along with them to obtain a product equals $1$ (depending on the parity of $n/2$), and then include pairs of elements of $S'$ just to complete $n/2$ terms. Consider the following cases and subcases:

\begin{enumerate}[(i)]
\item $n \ge 10$. We have $\frac{n}{2} - 3 \ge \frac{n-2}{4} \ge 2$. Consider the subcases according to $|S_1|$ and $|S_2|$.
\begin{enumerate}[({i.}1)]
\item 
Suppose that $|S_1| \ge 3$. Then at least one term from $S'$ comes from a product of two terms $(xy^{u_1}z) \bd (xy^{u_2}z) \mid S_1$ having the same parity. \\
We first consider that $(xy^{u_1}z) \cdot (xy^{u_2}z) = y^{2\alpha}$. If $\alpha \equiv 0 \pmod {n/2}$, then 
\begin{align*}
(xy^{u_1}z) \cdot (y^{2\beta}) \cdot (xy^{u_2}z) \cdot (y^{2\beta}) \cdot (y^{2\alpha})^{n/2-3} = 1 \quad &\text{ if $v_{y^{2\alpha}}(S') \ge \frac n2 - 2$, \; and } \\ 
(xy^{u_1}z) \cdot (y^{2\beta})^2 \cdot (xy^{u_2}z) \cdot (y^{2\beta})^2 \cdot (y^{2\alpha})^{n/2-5} = 1 \quad &\text{ if $v_{y^{2\alpha}}(S') = \frac n2-3$}.
\end{align*}
These products have $n$ terms (considering as terms of $S$). The case $\beta \equiv 0 \pmod {n/2}$ is similar, so we consider $\alpha, \beta \not\equiv 0 \pmod {n/2}$. \\
If $n/2$ is even, then 
$$(xy^{u_1}z) \cdot (y^{2\alpha}) \cdot (y^{2\beta})^{n/4-1} \cdot (xy^{u_2}z) \cdot (y^{2\beta})^{n/4-1} = 1$$ 
and this product has $n$ terms. \\
If $n/2$ is odd, then the set $\{\pm k(\alpha-\beta) \pmod {n/2} \mid 1 \le k \le \frac{n-2}{4} \}$ generates $n/2 - 1$ distinct and non-zero elements modulo $n/2$, since $\gcd(\alpha-\beta, n/2) = 1$. In particular, $\pm k(\alpha-\beta) \equiv -\alpha \pmod {n/2}$ for some $1 \le k \le \frac{n-2}{4}$, thus $\pm k(\alpha-\beta) - \alpha \equiv 0 \pmod {n/2}$ for some sign $\pm$. Thus either 
\begin{align*}
(xy^{u_1}z) \! \cdot \! (y^{2\beta})^{(n-2)/4} \! \cdot \! (xy^{u_2}z) \! \cdot \! (y^{2\alpha})^{k} \! \cdot \! (y^{2\beta})^{(n-2)/4 - k} &= 1 \\
\text{ for the positive sign } &+k, \; \text{ or} \\
(xy^{u_1}z) \! \cdot \! (y^{2\beta})^{(n-2)/4} \! \cdot \! (y^{2\alpha}) \! \cdot \! (xy^{u_2}z) \! \cdot \! (y^{2\alpha})^{k-1} \! \cdot \! (y^{2\beta})^{(n-2)/4 - k} &= 1 \\
\text{ for the negative sign } &-k,
\end{align*}
and these products have $n$ terms. \\
The subcases $(xy^{u_1}z) \cdot (xy^{u_2}z) = y^{2\beta}$ and $(xy^{u_1}z) \cdot (xy^{u_2}z) = y^{2\gamma}$ (where the latter only occurs if $|R| = 8$) are completely similar, as well as the case $|S_2| \ge 3$. \\

\item Suppose that $2 \le |S_1| + |S_2| \le 4$. 

If either $|S_1| = 2$ or $|S_2| = 2$ and the respective terms have the same parity, then we obtain one more term in $C_{n/2}$; in this case we are done if $|R| = 6$ since we would obtain $|S'| = n-1$ terms over $C_{n/2}$. Thus, suppose that $|R| = 8$. Then the product of this pair must be $y^{2\alpha}$ and we must have $S' = (y^{2\alpha})^{[n/2-1]} \bd (y^{2\beta})^{[n/2-1]}$ in order to avoid product-one of length $n/2$. If $2\alpha \not\equiv 0 \pmod {n/2}$, then we just change the order, obtaining another sequence which is not extremal. Otherwise, we may assume without loss of generality that $xy^{u_1}z \cdot xy^{u_2}z = y^{2\alpha} \in \{1,y^{n/2}\}$, where $y^{2\alpha} = y^{n/2}$ only if $n/2$ is even. Thus either 
\begin{align*}
\quad \quad \quad (xy^{u_1}z) \cdot (y^{2\beta})^{n/4-1} \cdot (xy^{u_2}z) \cdot (y^{2\beta})^{n/4-1} \cdot y^{2\alpha} &= 1 \quad \text{ if $n/2$ is even, \; or} \\
(xy^{u_1}z) \cdot (y^{2\beta})^{(n-2)/4} \cdot (xy^{u_2}z) \cdot (y^{2\beta})^{(n-2)/4} &= 1 \quad \text{ if $n/2$ is odd (that is, $y^{2\alpha} = 1$).}
\end{align*}

Therefore, we assume that no terms of $S'$ comes from $S_1 \bd S_2$, that is, the terms of $S_1$ and of $S_2$, whenever there are two, have opposite parities. Lemma \ref{termosR} ensures that there exists a subsequence $T = g_1 \bd g_2 \bd g_3 \bd g_4 \mid R$ with $g_1 \bd g_2 \mid S_1 \bd S_2$ such that $g_1 \cdot g_2 \cdot g_3 \cdot g_4 = y^{2\delta} \in C_{n/2}$. Moreover, either $2\delta \not\equiv 2\alpha \pmod n$ or $2\delta \not\equiv 2\beta \pmod n$. Suppose without loss of generality that $2\delta \not\equiv 2\alpha \pmod n$. Notice that $(y^{2\alpha})^{[n/2-3]} \bd (y^{2\beta})^{[n/2-2]} \mid S'$. The equation $\pm k(\alpha - \beta) \equiv \alpha-\delta \pmod {n/2}$ has a solution for some sign $\pm$ and some $1 \le k \le \frac{n-2}{4}$. It implies that 
\begin{align*}
\quad \quad \quad g_1 \cdot (y^{2\beta})^{(n-2)/4} \cdot g_2 \cdot g_3 \cdot g_4 \cdot (y^{2\alpha})^{k-1} \cdot (y^{2\beta})^{(n-2)/4 - k} &= 1 \\
\text{ for the positive sign } &+k, \; \text{ or} \\
g_1 \cdot (y^{2\alpha})^{k+1} \cdot (y^{2\beta})^{(n-6)/4 - k} \cdot g_2 \cdot g_3 \cdot g_4 \cdot (y^{2\beta})^{(n-6)/4} &= 1 \\
\text{ for the negative sign } &-k \text{ where $1 \le k < \frac{n-2}{4}$, \; or} \\
g_1 \cdot g_2 \cdot g_3 \cdot g_4 \cdot (y^{2\alpha})^{(n-6)/4} \cdot (y^{2\beta})^{(n-2)/4} &= 1 \\
\text{ for the negative sign } &-k \text{ where $k = \frac{n-2}{4}$,}
\end{align*}
and all these products have $n$ terms.

\item If $|S_1| + |S_2| \le 1$, then $|R| \le 5$, a contradiction.
\end{enumerate}

\item $n=8$. We have $|S| = 18$ and either $|R| = 6$ or $|R| = 8$. By Lemma \ref{termosR}, there exists $T \mid R \mid S$ with $|T| = 4$ yielding two distinct products in $C_4$. Notice that the subsets $\{y^{4\alpha}, y^{2\alpha+2\beta}, y^{4\beta}\}$ and $\{y^{2\alpha+2\beta}, y^{4\beta}, y^{2\beta+2\gamma}\}$ of $C_4$ both contain three distinct elements and can be obtained as a product of two terms of $S'$, since $\alpha \not\equiv \beta \pmod 2$. Therefore, the inverse of at least one of the elements of these subsets matches one of the products of $T$. It produces a product-one of length $8$.

\item $n=6$. We have $|S| = 14$ and either $|R| = 6$ or $|R| = 8$. 

If $R \nmid xy^{\alpha}z \bd xy^{\alpha+3}z \bd xy^{\alpha} \bd xy^{\alpha+3} \bd z \bd y^3z \bd 1 \bd y^3$, then Lemma \ref{termosR} ensures that there exists $T \mid R \mid S$ with $|T| = 4$ yielding two distinct products in $C_3$. By Pigeonhole Principle, one of these products belongs to either $\{y^{-2\alpha},y^{-2\beta}\} \subset C_3$ or $\{y^{-2\beta},y^{-2\gamma}\} \subset C_3$. Therefore, we multiply the terms of $T$ (in a suitable order) by the respective term of $S'$, obtaining $1$ as a product of $6$ terms of $S$.

If $R \mid xy^{\alpha}z \bd xy^{\alpha+3}z \bd xy^{\alpha} \bd xy^{\alpha+3} \bd z \bd y^3z \bd 1 \bd y^3$, then, by changing terms among $S'$ and $R$, we must have 
$$S = (xy^{\alpha}z)^{[a_1]} \bd (xy^{\alpha+3}z)^{[a_2]} \bd (xy^{\alpha})^{[a_3]} \bd (xy^{\alpha+3})^{[a_4]} \bd (z)^{[a_5]} \bd (y^3z)^{[a_6]} \bd (1)^{[a_7]} \bd (y^3)^{[a_8]},$$ 
where $a_1+ \dots + a_8 = 14$ (otherwise, by Lemma \ref{termosR}, we are done). By Pigeonhole Principle, $S$ contains three pairs of identical terms, and their respective product (in any order) is $1$.

\item $n=4$. We have $|S| = 10$ and either $|R| = 6$ or $|R| = 8$. By Lemma \ref{termosR}, there exists $T \mid R \mid S$ with $|T| = 4$ yielding two distinct products in $C_2$. One of these products is $1$, thus we are done.
\end{enumerate}
It completes the proof.

\qed

\begin{theorem}
Let $n \ge 4$ be an even integer. It holds that $\s(D_{2n} \times C_2) = 2n+2$ and $\e(D_{2n} \times C_2) = n+3$. In particular, $\s(D_{2n} \times C_2) = \e(D_{2n} \times C_2) + \exp(D_{2n} \times C_2) - 1$.
\end{theorem}

\proof
By Lemma \ref{lowerbound}{\it (b)}, Proposition \ref{sle}, and Ineq. \eqref{sexp}, it follows that
$$n+3 \le \e(D_{2n} \times C_2) \le \s(D_{2n} \times C_2) - \exp(D_{2n} \times C_2) + 1 \le 2n+2 - n + 1 = n+3,$$
therefore the equalities hold.

\qed

\section{Proof of Conjecture $\E(G) = \d(G) + |G|$ for $G \simeq D_{2n} \times C_2$}\label{conjEdord}

\begin{proposition}\label{Ele}
Let $n \ge 4$ be an even integer. Then $\E(D_{2n} \times C_2) \le 5n+1$.
\end{proposition}

\proof
The argument follows the same arguments as in Proposition \ref{sle}. Let $S \in \F(D_{2n} \times C_2)$, with $|S| = 5n+1$. If $|S_1| + |S_2| \le 1$, then $|S_3| + |S_4| \ge 5n = 2n + \E(C_n \times C_2)$, therefore there exists $T_1 \bd T_2 \mid S$ with $|T_1| = |T_2| = 2n$ such that $1 \in \pi(T_1 \bd T_2)$. From now on, suppose that $|S_1| + |S_2| \ge 2$. 

Let $g' \bd g'' \mid S_1 \bd S_2$, and consider $S_0 = S \bd (g' \bd g'')^{[-1]}$. According to our notation, let $S_0 = S_{01} \bd S_{02} \bd S_3 \bd S_4$, where $S_{01} \mid S_1$, $S_{02} \mid S_2$, and $S_{01} \bd S_{02} = (S_1 \bd S_2) \bd (g' \bd g'')^{[-1]}$. As in the proof of Proposition \ref{sle}, it is possible to find 
$$\left\lceil \frac{|S_{01}| - 2}{2} \right\rceil + \left\lceil \frac{|S_{02}| - 2}{2} \right\rceil + \left\lceil \frac{|S_3| - 2}{2} \right\rceil + \left\lceil \frac{|S_4| - 2}{2} \right\rceil \ge \frac{|S_0| - 7}{2} = 5 \cdot \frac n2 - 4$$
elements in $C_{n/2}$, each of them being a product of exactly two terms of $S$. 

\begin{enumerate}[(i)]
\item {\bf Case} $n \ge 6$. Since 
\begin{equation}\label{ineqj}
5 \cdot \frac n2 - 4 - j \cdot \frac n2 \ge 2 \cdot \frac n2 - 1 
\end{equation}
for $j \in \{0,1,2\}$, the Erd\H os-Ginzburg-Ziv Theorem (first part of Lemma \ref{inverseegz}) implies that there exists a product-one subsequence of length $3 \cdot n/2$ over $C_{n/2}$, which corresponds to $3n$ terms over $S_0$. Let $S' \mid S_0$ be the subsequence formed by these $3n$ terms such that $1 \in \pi(S')$. In addition, let $S'' = S \bd S'^{[-1]} = (g' \bd g'' \bd S_0) \bd S'^{[-1]}$ so that $|S''| = 2n+1$, where $S'' = S''_1 \bd S''_2 \bd S''_3 \bd S''_4$, $S''_i \mid S_i$ for $i \in \{1,2,3,4\}$, and $|S''_1| + |S''_2| \ge 2$. Again as in the proof of Proposition \ref{sle}, we obtain at least
$$\left\lceil \frac{|S''_1| - 2}{2} \right\rceil + \left\lceil \frac{|S''_2| - 2}{2} \right\rceil + \left\lceil \frac{|S''_3| - 2}{2} \right\rceil + \left\lceil \frac{|S''_4| - 2}{2} \right\rceil \ge \frac{|S''| - 7}{2} = 2 \cdot \frac n2 - 3$$
elements in $C_{n/2}$, each of them being a product of exactly two terms of $S''$. Let $R \mid S'' \mid S$ be the sequence formed by the terms of $S''$ which are not used to form the (at least) $n - 3$ pairs, so that $|R| \in \{1,3,5,7\}$. Looking more closely, we notice that the exact number of pairs as the latter is $\frac{2n+1 - |R|}{2} = n - \frac{|R| - 1}{2}$. 

If $|R| \le 3$, then we can use the Erd\H os-Ginzburg-Ziv Theorem once again, which implies that there exists a subsequence of $S$ of length $4n$ whose product is $1$. 

If $|R| = 7$, then we may suppose that the last $n-3$ pairs over $C_{n/2}$ must comprise either 
\begin{align*}
S''' &= (y^{2\alpha})^{[n/2-2]} \bd (y^{2\beta})^{[n/2-1]}, \\ 
S''' &= (y^{2\alpha})^{[n/2-2]} \bd (y^{2\beta})^{[n/2-2]} \bd (y^{2\gamma}), \quad \text{or} \\ 
S''' &= (y^{2\alpha})^{[n/2-3]} \bd (y^{2\beta})^{[n/2-1]} \bd (y^{2\gamma}),
\end{align*}
where $\gcd(\alpha-\beta, n/2) = 1$. In this case, we will follow the steps of Proposition \ref{sle} (with $|R| = 8$) work, namely Cases (i.1), (i.2), (ii), and (iii). 

If $|R| = 5$, the last $n-2$ pairs must comprise the sequence $S''' = (y^{2\alpha})^{[n/2-1]} \bd (y^{2\beta})^{[n/2-1]} \in C_{n/2}$, where $\gcd(\alpha-\beta, n/2) = 1$. There are two subcases:

\begin{enumerate}[{(i.}1{)}]
\item {\bf Subcase:} $R$ contains exactly one term from $S_1 \bd S_2$. It implies that at least one term of $S'''$ is a product of a pair of $S_1 \bd S_2$, say $xy^{u_1}z \cdot xy^{u_2}z = y^{2\alpha}$.

If $n \ge 10$, then the proposition follows from the same argument than Cases (i.1) and (i.2 - ``the respective terms have the same parity'') of Proposition \ref{sle}.

If $n = 8$, then we consider two more cases. If $\alpha$ is odd, then we just change the order, obtaining a sequence of maximal length which is distinct from $(y^{2\alpha})^{[3]} \bd (y^{2\beta})^{[3]}$; therefore $\alpha$ is even. It implies that $xy^{u_1}z \cdot y^{2\beta} \cdot xy^{u_2}z \cdot y^{2\beta} \cdot y^{2\alpha} = 1$, and this product has $8$ terms from $S$. Thus we are done.

If $n = 6$, then the same argument implies that we must have $\alpha \equiv 0 \pmod 3$, that is, $xy^{u_1}z \cdot xy^{u_2}z = 1$ (otherwise, we would change the order). Since $xy^{u_1}z \cdot y^{2\beta} \cdot xy^{u_2}z \cdot y^{2\beta} = 1$ and this product has $6$ terms from $S$, we are done.

\item {\bf Subcase:} $R$ contains at least two terms from $S_1 \bd S_2$. In particular, we are in the hypothesis of Lemma \ref{termosR} (b). From the previous subcase, we may assume that every term of $S'''$ comes from pairs of $S_3 \bd S_4$.

If $n \ge 10$, then the proposition follows from the same arguments than Case (i.2 - ``no terms of $S'$ comes from $S_1 \bd S_2$'') of Proposition \ref{sle} (with $|R| = 6$), using Case (b) of Lemma \ref{termosR}.

If $n = 8$, then our argument follows the same steps than Case (ii) of Proposition \ref{sle} (with $|R| = 6$), using Case (b) of Lemma \ref{termosR}.

If $n = 6$, then Lemma \ref{termosR} implies that there exists $T \mid R$ with $|T| = 4$ such that $\pi(T) \subset C_{n/2}$. If $T \nmid xy^{\alpha}z \bd xy^{\alpha+3}z \bd xy^{\alpha} \bd xy^{\alpha+3} \bd z \bd y^3z \bd 1 \bd y^3$, then we use the same argument than Case (iii) of Proposition \ref{sle}, thus we are done. On the other hand, if $T \mid xy^{\alpha}z \bd xy^{\alpha+3}z \bd xy^{\alpha} \bd xy^{\alpha+3} \bd z \bd y^3z \bd 1 \bd y^3$, then $T$ has two or four terms of $S_1 \bd S_2$ and $\pi(T)$ has only one term. Furthermore, since every term of $S'''$ comes from pairs of $S_3 \bd S_4$, it follows that $9 \le |S_3| + |S_4| \le 11$. Suppose without loss of generality that $y^{w_1}z \cdot y^{w_2}z = y^{2\alpha}$ (the case $y^{t_1} \cdot y^{t_2} = y^{2\alpha}$ is completely similar). In any case, we obtain $\pi(T) = \{1\}$, therefore we may assume that $0 \not\in \{2\alpha,2\beta\}$ modulo $6$ (otherwise we are done). It implies that $y^{w_1 \pm w_2} \in \pi(T \bd y^{w_1}z \bd y^{w_2}z)$, therefore the terms of $S_3$ are all distinct, as well as the terms of $S_4$. Since $|S_3| + |S_4| \ge 9$, it follows that either $|S_3| \ge 5$ or $|S_4| \ge 5$. Suppose without loss of generality that $|S_3| \ge 5$. Thus we can set either $\{w_1,w_2\} = \{1,5\}$ modulo $n$ or $\{w_1,w_2\} = \{2,4\}$ modulo $n$, which implies that $0 \equiv w_1+w_2 \equiv 2\alpha \pmod n$, a contradiction. Hence we are done.
\end{enumerate}

\item {\bf Case} $n = 4$. Notice that $|S| = 21$, $C_2 \simeq \{1,y^2\}$, and $\s(C_2) = 3$. Since 
\begin{equation}\label{sn4}
\left\lceil \frac{|S_1| - 2}{2} \right\rceil + \left\lceil \frac{|S_2| - 2}{2} \right\rceil + \left\lceil \frac{|S_3| - 2}{2} \right\rceil + \left\lceil \frac{|S_4| - 2}{2} \right\rceil \ge \frac{21 - 7}{2} = 7 \ge 3 + 2j
\end{equation}
for $j \in \{0,1,2\}$, the Erd\H os-Ginzburg-Ziv Theorem (first part of Lemma \ref{inverseegz}) implies that there exist three disjoint subsequences of $S$, say $S'$, $S''$ and $S'''$, such that $1 \in \pi(S') \cap \pi(S'') \cap \pi(S''')$ and $|S'| = |S''| = |S'''| = 4$. It remains at least one pair of terms over $S_0 = S \bd (S' \bd S'' \bd S''')^{[-1]}$ whose product belongs to $C_2$. Let $R \mid S \bd (S' \bd S'' \bd S''')^{[-1]}$ be the sequence of unpaired terms, so that $|R| \in \{1,3,5,7\}$. 

If $|R| \le 3$, then the left hand side of Ineq. \eqref{sn4} is at least $\frac{21 - 3}{2} = 9 \ge 3 + 2j$ for $j \in \{0,1,2,3\}$. Thus there exists $S'''' \mid S_0$ with $|S''''| = 4$ such that $1 \in \pi(S'''')$, therefore $1 \in \pi(S' \bd S'' \bd S''' \bd S'''')$. 

If $|R| = 7$, then, by Lemma \ref{termosR}, there exists $T \mid R$ with $|T| = 4$ such that $\pi(T)$ contains at least two distinct elements of $C_2$. In particular, $1 \in \pi(T)$, hence $1 \in \pi(S' \bd S'' \bd S''' \bd T)$.

If $|R| = 5$, then it remains two pairs over $S_0$ whose product belongs to $C_2$. These pairs, namely $g' \bd g''$ and $g''' \bd g''''$, must comprise the products (say) $g' \cdot g'' = 1$ and $g''' \cdot g'''' = y^2$ (otherwise, the pairs form equal products over $C_2$, hence we obtain another sequence $S''''$ of length $4$ such that $1 \in \pi(S'''')$, and we are done). Furthermore, by Lemma \ref{termosR}, there exists $T \mid R$ with $|T| = 4$ such that the product of the terms of $T$ belongs to $C_2$. If this product is $1$, then $1 \in \pi(S' \bd S'' \bd S''' \bd T)$, thus we are done. If this product is $y^2$, then $1 \in \pi(S' \bd S'' \bd g' \bd g'' \bd T \bd g''' \bd g'''')$. The latter product has $16$ terms, thus we are done.
\end{enumerate}
It completes the proof.

\qed

\begin{theorem}
Let $n \ge 4$ be an even integer. It holds that $\E(D_{2n} \times C_2) = 5n+1$ and $\d(D_{2n} \times C_2) = n+1$. In particular, $\E(D_{2n} \times C_2) = \d(D_{2n} \times C_2) + |D_{2n} \times C_2|$.
\end{theorem}

\proof
By Lemma \ref{lowerbound}{\it (a)}, Proposition \ref{Ele}, and Ineq. \eqref{Ed}, it follows that
$$n+1 \le \d(D_{2n} \times C_2) \le \E(D_{2n} \times C_2) - |G| \le 5n+1 - 4n = n+1,$$
therefore the equalities hold.

\qed

\end{document}